\documentclass[aps,prl,reprint,amsmath,amsfonts,amssymb, superscriptaddress]{revtex4-1}
\usepackage{amsthm}
\usepackage{amsmath}
\usepackage{graphicx}
\usepackage{color}
\usepackage{enumitem}
\usepackage[hidelinks]{hyperref}
\usepackage{xcolor}
\hypersetup{
    colorlinks,
    linkcolor={blue!50!black},
    citecolor={blue!50!black},
    urlcolor={blue!80!black}
}

\renewcommand{\vec}[1]{\mathbf{#1}}

\newtheorem{theorem}{Theorem}

\begin{document}
\title{Chebyshev approximation and the global geometry of sloppy models}

\author{Katherine N. Quinn}
\affiliation{Physics Department, Cornell University, Ithaca, NY 14853-2501, United States. (knq2@cornell.edu)}

\author{Heather Wilber}
\affiliation{Center for Applied Mathematics, Cornell University, Ithaca, NY 14853-3801, United States. (hdw27@cornell.edu)}

\author{Alex Townsend}
\affiliation{Mathematics Department, Cornell University, Ithaca, NY 14853-4201, United States. (townsend@cornell.edu)}

\author{James P. Sethna}
\affiliation{Physics Department, Cornell University, Ithaca, NY 14853-2501, United States. (sethna@lassp.cornell.edu)}

\date{\today}

\begin{abstract}
Sloppy models are complex nonlinear models with outcomes that are significantly affected by only a small subset of parameter combinations. Despite forming an important universality class and arising frequently in practice,  formal and systematic explanations of sloppiness are lacking. 
By unifying geometric interpretations of sloppiness with Chebyshev approximation theory, we offer such an explanation, and show how sloppiness can be described explicitly in terms of model smoothness. Our approach results in universal bounds on model predictions for classes of smooth models, and our bounds capture global geometric features that are intrinsic to their model manifolds.  We illustrate these ideas using three disparate models: exponential decay, reaction rates from an enzyme-catalysed chemical reaction, and an epidemiology model of an infected population.
\end{abstract}

\maketitle

Complex nonlinear models used to simulate and predict experimentally observed phenomena often exhibit a structural hierarchy: Perturbing some model parameter combinations drastically impacts predictions, whereas others can vary widely without effect. Such models are called \textit{sloppy}. Sloppy models appear to be common, arising in systems biology~{\cite{Transtrum2003,Brown2004, Gutenkunst_SloppyParameters_2007}}, accelerator physics~\cite{Ryan2007}, radioactive decay~\cite{Ruhe1980}, critical phenomena~\cite{Machta_ParameterCompression_2013}, insect flight~\cite{Berman2007}, and many other areas~\cite{Transtrum2015}. 
Developing a formal and rigorous framework to explain sloppiness is essential since it leads to systematic methods of low-dimensional approximation and parameter fitting. In this letter, we unify recently developed geometric descriptions of sloppiness~\cite{Transtrum2015} with classical ideas from polynomial approximation theory~\cite{trefethen2013approximation}. We posit that in many cases, sloppiness is fundamentally linked to the smoothness of the model, and provide a rigorous description of this connection. 

Sloppy models are characterized by the geometry of their \textit{model manifolds}, i.e., the space of all possible predictions. Studying this geometry yields fruitful information for several reasons: (1) the dominant components reflect emergent behavior of the models (how the microscopic interactions do or do not produce macroscopic behavior~\cite{Machta_ParameterCompression_2013}), (2) the boundaries represent reduced-model approximations~\cite{Transtrum2014}, and (3) knowledge of the manifold geometry leads to more efficient data fitting methods~\cite{Transtrum2011}. Model manifolds take on the shape of \textit{hyperribbons}~\cite{Transtrum2010,Transtrum2011}, so-called because they resemble high-dimensional ribbons. They are much longer than they are wide, much wider than they are thick, etc., yielding effective low-dimensional representations ~\cite{Transtrum2010,Transtrum2011,Transtrum2015}.  In~\cite{Transtrum2010}, polynomial interpolation is used to show that certain analytic models must have manifolds that are hyperribbons. Drawing on this work, we apply more general methods of polynomial approximation and develop two key results:  (1) we derive explicit geometric bounds that explain why many multiparameter model manifolds are hyperribbons, and (2) we numerically bound the prediction space of all models that share certain features of smoothness. 

Consider a nonlinear model that depends continuously on input parameters \smash{$\theta= (\theta^1,\dots,\theta^K)$}  to generate predictions \smash{$y_\theta(t)$}. We assume without loss of generality that the model is shifted and scaled so that $t \in  [-1, 1]$.   
Using model predictions \smash{$y_\theta(t_k) = Y_k$} at points \smash{$\{t_k\}_{k= 0}^{N-1} \subset [-1, 1]$}, $N \geq K$, the model manifold $\mathcal{Y}$ is defined as the $K$-dimensional surface parameterized by \smash{$Y(\theta) =  (Y_0, \ldots, Y_{N-1})$} and embedded in the $N$-dimensional prediction space.  

To bound the manifold $\mathcal{Y}$ and study its geometry, we study manifolds associated with polynomial approximations to \smash{$y_\theta$}.  Let \smash{$\{\phi_j\}_{j = 0}^{\infty}$} be a complete polynomial basis, and suppose that \smash{$y_\theta(t)  =  \sum_{j = 0}^{\infty} b_j(\theta) \phi_j(t)$.}   Let \smash{$P(\theta)  = (P_0, \ldots, P_{N-1})$} define the model manifold $\mathcal{P}$ of \smash{$p_{N-1}(\vec{t}; \theta) = \sum_{j = 0}^{N-1} b_j(\theta) \phi_j(\vec{t})$}, where \smash{$\vec{t} = (t_0, \ldots, t_{N-1})^T$.}  We view the coefficients \smash{$(b_0(\theta), \ldots, b_{N-1}(\theta))$} as a set of $N$ parameters, with model outcomes given by \smash{$P_k = p_{N-1}(t_k; \theta)$}.  

By definition, $P(\theta) = X\vec{b}$, where \smash{$X_{ij} = \phi_{j-1}(t_{i-1})$} and \smash{$\vec{b} = (b_0, \ldots, b_{N-1})^T$}. Here, $X$ forms a linear map from parameter space to prediction space. We characterize the geometry of $\mathcal{P}$ using the singular values of $X$.  Suppose, for example,  that \smash{$\|\vec{b}\|_2 < r$}, so that the parameter space is bounded in $S$, an $N$-sphere of radius $r$. The action of $X$ on $S$ distorts it into a hyperellipsoid \smash{$H_P$}. If \smash{$\ell_j(H_P)$} is the diameter of the $j$th cross-section of \smash{$H_P$},  then 
\begin{equation} 
\label{eq:dim_bound}
\ell_j(H_P)  = 2 r \sigma_j(X),
\end{equation}
where \smash{$\sigma_1(X) \geq \sigma_2(X) \geq \cdots \geq \sigma_{N}(X)$} are the singular values of $X$.  
When $X$ has rapidly decaying singular values, \smash{$H_P$} has a hyperribbon structure. Accounting for the error $\|y_\theta - p_{N-1}\|_\infty$,  where $\|\cdot\|_{\infty}$ is the $L^\infty$ norm on $[-1, 1]$, we find that a hyperellipsoid $H_Y$ must enclose $\mathcal{Y}$, with cross-sectional widths given by
\begin{equation} 
\label{eq:dim_bound2}
\ell_j(H_Y)  = 2 r \sigma_j(X) + 2\|y - p_N\|_\infty.
\end{equation}
We consider two applications of this idea: First, we choose our basis functions \smash{$\{\phi_j\}_{j = 0}^\infty$} as the Chebyshev polynomials. 
Truncated Chebyshev expansions converge to \smash{$y_\theta$} at an asymptotically optimal rate for polynomial approximation~\cite{trefethen2013approximation}. As we show below, this rate controls the magnitude of \smash{$\sigma_j(X)$} in Eq.~\eqref{eq:dim_bound2}, and can be used to explicitly bound the cross-sectional widths of \smash{$H_Y$}. 
We also analyze the case where $\{\phi_j\}_{j = 0}^\infty$ are the monomials and \smash{$p_{N-1}$} is the truncated Taylor series expansion of \smash{$y_\theta$}. In this case, we observe that the numerical computation of \smash{$\sigma_j(X)$} results in excellent practical and universal bounds on the prediction space for large classes of models. 

\vspace{.25cm}

 \textbf{Chebyshev expansions.} Suppose that \smash{$y_\theta$} has a convergent Chebyshev expansion, so that it is given by \smash{$y_\theta(t) = \sum_{j = 0}^\infty	c_j(\theta) T_j(t)$}, where \smash{$T_j(t) = \cos(j \arccos t)$} is the degree $j$ Chebyshev polynomial~\cite[Ch.~3]{trefethen2013approximation}.   We can approximate \smash{$y_\theta$} with a degree $\leq N-1$ polynomial by truncating the Chebyshev series after $N$ terms:
\begin{equation}
\label{eq:cheb_appx}
p_{N-1}(t; \theta) = \sum_{j = 0}^{N-1} c_j(\theta) T_j(t).
\end{equation}
Truncated Chebyshev expansions have near-best global approximation properties. The error \smash{$\|y_\theta - p_{N-1}\|_\infty$} is within a $\log N$ factor of \smash{$\|y_\theta - p^{best}_{N-1}\|_\infty$}~\cite[Ch.~16]{trefethen2013approximation}, where \smash{$p_{N-1}^{\text{best}}$} is the best polynomial approximant to \smash{$y_\theta$} of degree $\leq N\!- \! 1$. We cannot directly use \smash{$p^{best}_{N-1}$} in our arguments because bounds on \smash{$\|y_\theta - p^{best}_{N-1}\|_\infty$} are only known in an asymptotic sense. Fortunately, explicit bounds on \smash{$\|y_\theta - p_{N-1}\|_\infty$} are known when \smash{$y_\theta$} is sufficiently smooth.

 We first consider the case where \smash{$y_\theta$} is analytic in an open neighborhood of $[-1, 1]$. Such a region contains a \textit{Bernstein ellipse} $E_\rho$, defined as the image of the circle $|z| = \rho$ under the Joukowsky mapping \smash{$(z +z^{-1})/2$}. It has foci at $\pm 1$, and the lengths of its semi-major and semi-minor axes sum to $\rho$. The polynomial in Eq.~\eqref{eq:cheb_appx} converges to \smash{$y_\theta$} as $N \to \infty$  at a  rate determined by $\rho$: 
\begin{theorem} 
\label{thm:Analytic_conv}
Let $M> 0$ and $\rho > 1$ be constants and suppose that \smash{$y_\theta(t)$}, $t \in [-1, 1]$, is analytically continuable to the region enclosed by the Bernstein ellipse \smash{$E_\rho$},  with \smash{$|y_\theta | \leq M$} in \smash{$E_\rho$},  uniformly in $\theta$.  Let \smash{$p_{N-1}(t ; \theta)$} be as  in Eq.~\eqref{eq:cheb_appx}. Then, 
\begin{align}
\label{eq:analyticbound}
&  (i) \, \, \left \| y_\theta - p_{N-1} \right \|_{\infty} \leq \frac{2M \rho^{-N+1}}{\rho-1}, \\
\label{eq: coeff_bounds}
 &  (ii) \, \, |c_0| \leq M, \, \, |c_j(\theta)| \leq 2M\rho^{-j}, \qquad j \geq 1. 
 \end{align}
\end{theorem}
\begin{proof} For a proof, see Theorem 8.2 in~\cite{trefethen2013approximation}. \end{proof}

To exploit the decay of the coefficients in Eq.~\eqref{eq: coeff_bounds}, we define  modified coefficients \smash{$\tilde{c}_j = \rho^j c_j$}. We then have that \smash{$P(\theta) =  X \vec{\tilde{c}}$}, where $X = JD$,  \smash{$J_{ij} = T_{j-1}(t_{i-1})$}, and $D$ is diagonal with entries \smash{$D_{jj} = \rho^{-(j-1)}$}. By~\eqref{eq: coeff_bounds}, we have that \smash{$\|\vec{\tilde{c}}\|_2 < 4M\sqrt{4N-3}$}. This implies that $\mathcal{P}$ is bound in a hyperellipsoid \smash{$H_P$}. By Eq~\eqref{eq:dim_bound}, we have that \smash{$\ell_j(H_P) = 8 M \sqrt{4N-3} \sigma_j(X)$}. 
To bound \smash{$\sigma_j(X)$} explicitly, we first prove a conjecture proposed in~\cite{Waterfall2006}:
\begin{theorem} 
\label{thm: ESE}
Let \smash{$S\in\mathbb{R}^{N\times N}$} be symmetric and positive definite. Let \smash{$E\in\mathbb{R}^{N \times N}$} be diagonal with \smash{$E_{ii} = \epsilon^{i-1}$} and $0<\epsilon<1$. If \smash{$\lambda_1 \geq \lambda_2 \geq\cdots \geq \lambda_{N}$} are the ordered eigenvalues of $ESE$, then \smash{$\lambda_{m+1} = \mathcal{O}(\epsilon^{2m})$}. Specifically, 
\begin{equation} 
\lambda_{m+1} \leq \frac{\epsilon^{2m}}{1 - \epsilon^2} \max_{1\leq j,k\leq N} \left|S_{jk}\right|, \qquad 1 \leq m \leq N-1. 
\end{equation}  
\end{theorem} 
\begin{proof}\footnote{Previous proofs with weaker bounds were provided through private communications with Ari Turner and Yaming Yu.}
Consider the rank $m$ matrix
\begin{equation}
S_m = S(:,1\!:\!m) S(1\!:\!m,1\!:\!m)^{-1}S(1\!:\!m,:),
\end{equation} 
where $1 \!\leq \!m \!\leq N\!-\!1,$ and the notation $M(\,: \, , 1\!:\!m)$ denotes the submatrix of $M$ consisting of its first $m$ columns. 
Clearly, \smash{$S_m$} is well-defined because ${S(1\!:\!m,1\!:\!m)}$ is a principal minor of a positive definite matrix and is therefore invertible. Moreover, it can be verified that $(S - S_m)_{jk}  = 0$ for $1\leq j,k\leq m$.

Since $ESE$ is positive definite and \smash{${\rm rank}(S_m)= m$}, we know that \smash{$\lambda_{m+1} \leq \|E(S-S_m)E\|_2$}, where \smash{$\|\cdot\|_2$} denotes the spectral matrix norm~\cite[Ch.~2]{Gene1996}. Using \smash{$\|\cdot\|_F$} to denote the Frobenius norm, we have
\[
\begin{aligned}
\lambda_{m+1}^2 &\leq \|E(S-S_m)E\|_2^2 \leq \|E(S-S_m)E\|_F^2 \\
&= \sum_{j=m+1}^{N}\sum_{k=m+1}^{N} \epsilon^{2(j-1)+2(k-1)} \left|S_{jk} - (S_m)_{jk}\right|^2 \\
&\leq \frac{\epsilon^{4m}}{(1-\epsilon^2)^2}\max_{1\leq j,k\leq N} \left|S_{jk} - (S_m)_{jk}\right|^2 \\
&\leq \frac{\epsilon^{4m}}{(1-\epsilon^2)^2}\max_{1\leq j,k\leq N} \left|S_{jk}\right|^2 
\end{aligned} 
\]  
where the last inequality comes from the fact that the block \smash{$S(m\!+\!1\!:\!N, m\!+\!1\!:\!N)-S_m(m\!+\!1\!:\!N, m\!+\!1\!:\!N)$} is the Schur complement of $S(1\!:\!m, 1\!:\!m)$ in $S$~\cite{Gene1996}. 
\end{proof} 
Applying Theorem 2 to  \smash{$X^TX = DJ^TJD$}, we have that for $j > 1$,  \smash{$\sigma_j(X) \leq \sqrt{N} \rho^{-j+2}/\sqrt{\rho^2-1}$}, 
where we have used the fact that \smash{$|T_k(t)|\leq 1$} for $k\geq 0$ and $-1\leq t\leq 1$. It follows from Equations~\eqref{eq:dim_bound2} and~\eqref{eq:analyticbound} that predictions for \smash{$y_\theta(t)$} are  bounded by a hyperellipsoid \smash{$H_Y$}, with 
\begin{eqnarray}
\label{eq:HY}
\ell_j(H_Y) \leq \frac{2M\sqrt{4N^2-3N}\rho^{-j+2}}{\sqrt{\rho^2 - 1}}  +  \frac{4M \rho^{-N+1}}{\rho-1}, 
\end{eqnarray}
for $2\leq j \leq N.$
These bounds indicate that the hyperribbon structure of \smash{$H_Y$} is controlled by $\rho$, a parameter characterizing the analyticity of the model.  As $\rho$ becomes larger, bounds on the widths of the successive cross-sections of \smash{$H_Y$}  must decay more rapidly: In principle, \smash{$H_Y$} becomes successively thinner and more  ribbon-like.  

When \smash{$y_\theta$} is not analytic on an open neighborhood of $[-1, 1]$,  the decay rate of $\sigma_j(JD)$ is instead controlled by the smoothness of \smash{$y_\theta$} on $[-1, 1]$. We provide more discussion in the supplementary materials.

\vspace{.25cm}

\textbf{Taylor expansions.} The degree $N\! -\!1$ truncated Taylor polynomial of \smash{$y_\theta$} is  \smash{$p_{N-1}(t) = \sum_{k = 0}^{N-1} a_k(\theta)  (t-t_0)^j$}, where \smash{$a_k(\theta)  = y_\theta^{(k)}(t_0)/k!$}. 
We describe the analyticity of \smash{$y_\theta$} using the following condition: For all $N \geq 1$, 
\begin{eqnarray}
\label{eq:deriv_bnd}
\sum_{k=0}^{N-1} \left(\frac{R^k}{k!}\frac{d^ky_\theta(t)}{dt^k}\right)^2 < C^2 N
\end{eqnarray}
where $C > 0$,$R > 1$ are constants in $\theta$.  It follows that
\smash{$\|  d^{k}y_\theta /dt^k /k! \|_\infty < C\sqrt{k+1}R^{-k}$} for $k \geq 0$. The Taylor series for \smash{$y_\theta$} expanded about any point $t \in [-1, 1]$ has a radius of convergence of at least $R$. 
If $t_0 = 0$ and $R > 1$, then we find by simple estimates that
\begin{equation} 
\label{eq:taylor_error}
\|y - p_{N-1}\|_\infty \leq  \frac{C(NR -N + R)}{(1-R)^2} R^{-N+1}.
\end{equation} 
 As with the Chebyshev coefficients, we define \smash{$\tilde{a}_k = R^{k} a_k$}, and write \smash{$P(\theta) = V D \vec{\tilde{a}}$}, where \smash{$V_{ij} = t_{i-1}^{j-1}$} and \smash{$D = {\rm diag}(R^0, \ldots,  R^{-(N-1)})$}.  Explicit bounds on the singular values of $VD$ can be derived using its displacement structure~\cite{beckermann2017singular}. However, we require bounds that are characterized by the analyticity of \smash{$y_\theta$}. For this reason, we instead  apply Theorem~\ref{thm: ESE} to \smash{$DV^TVD$}, so that $\sigma_j(VD)$ is bounded in terms of $R$. By applying the constraint from Eq.~\eqref{eq:deriv_bnd} to $p_{N-1}$, we see that \smash{$\|\vec{\tilde{a}}\|_2 < C \sqrt{N}$}. It follows that the manifold $\mathcal{P}$ is bounded in a hyperellipsoid \smash{$H_P$}, where for $j \geq 2$,
 \begin{equation} 
 \label{eq:vand_bnds}
 \ell_j(H_P) = 2C \sqrt{N} \sigma_{j}(VD) \leq \frac{2CN}{ \sqrt{R^2 - 1}} R^{-j + 2} .
 \end{equation}
One can then form \smash{$H_Y$} in Eq.~\eqref{eq:dim_bound2}  by combining Eq.~\eqref{eq:taylor_error} and Eq.~\eqref{eq:vand_bnds}. 
How do these bounds compare to the Chebyshev-based results? The constraint in Eq.~\eqref{eq:deriv_bnd} implies that \smash{$y_\theta$} is analytic in the region $\mathcal{R}$ of the complex plane of distance $< R$ from $[-1, 1]$. 
It can be shown that \smash{$y_\theta$} must also be analytic and bounded by a function \smash{$M(\zeta)$} on any Bernstein ellipse \smash{$E_{\rho(\zeta)}$} in $\mathcal{R}$, with  \smash{$\rho(\zeta) = \zeta+ \sqrt{\zeta^2 + 1}$}~\cite{demanet2016stable}. The largest such ellipse is given by \smash{$\rho_{\max} = R + \sqrt{R^2 +1}$}, suggesting that Chebyshev-based bounds can  improve~\eqref{eq:vand_bnds} by nearly a factor of \smash{$2^j$}. However, $M(\zeta)$  is unbounded as $\zeta \to R$, so one must select $0 < \zeta< R$ to minimize the Chebyshev bound. Even when $\zeta$ is selected carefully, the conversion from Eq.~\eqref{eq:deriv_bnd} to a constraint involving \smash{$E_{\rho(\zeta)}$} may introduce an unphysically large constant into the bound.

The above argument implies that the \textit{explicit} bound from Eq.~\eqref{eq:vand_bnds} has a decay rate that is suboptimal for describing the hyperribbon widths of $\mathcal{Y}$.  However, we take advantage of the following fascinating observation: For moderate $N$, the numerically computed singular values \smash{$\sigma_{j}(VD)$} often decay more rapidly than the bound in Eq.~\eqref{eq:vand_bnds} captures. In fact, numerical tests show that their decay involves a kink: For small to moderate $j$, the magnitude of \smash{$\sigma_{j}(VD)$} often appears to decay at the Chebyshev-based rate \smash{$\mathcal{O}(\rho_{\max}^{-j})$}. When $j$ is larger, \smash{$\sigma_{j}(VD)$} instead decays at the weaker rate \smash{$\mathcal{O}(R^{-j})$}.  The supplementary materials include an expanded discussion of this phenomenon. 
Due to this behavior, the hyperellipsoid \smash{$H_Y$} formed by using the singular values of $VD$ directly results in a practical, universal boundary for the prediction space of all models obeying Eq.~\eqref{eq:deriv_bnd}.

\begin{figure*}
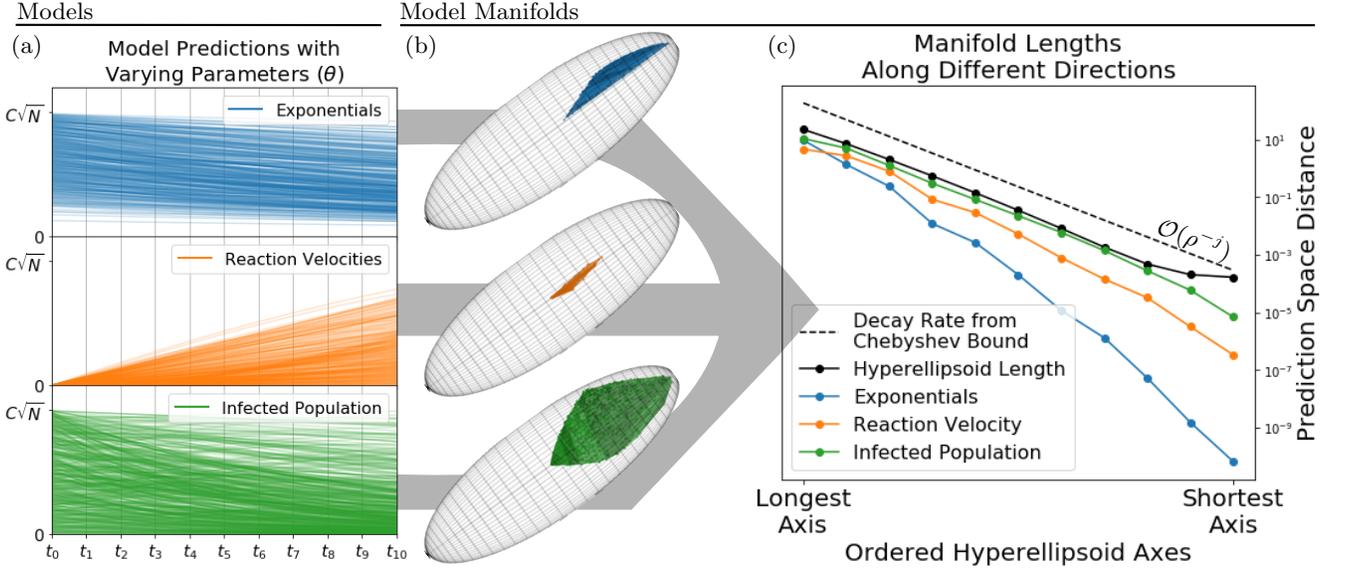
   \caption{\textbf{Model manifold} of three disparate models: (1) exponential decay, (2) reaction velocities of an enzyme-catalysed reaction, and (3) the infected in an SIR model. The models are evaluated at $11$ equally spaced points on $[0, 1]$ (shifted and scaled from $[-1, 1]$), and obey the smoothness condition in Eq.~\eqref{eq:deriv_bnd}, with  $C = 1$ and $R = 2$. (a)~An illustration of each model, where each line represents the respective model predictions with a different set of parameters. (b)~Each point in the figure displays the model prediction for a given parameter choice. The space of all possible predictions forms a geometric object known as the \textit{model manifold} (shown as the colored shape in the middle of each ellipse). The model manifolds are all bounded by the same hyperellipsoid, and so the two axes represent the first and second longest hyperellipsoid axes. Note that, in all three models, only values greater than 0 are physically significant. This constraint manifests itself geometrically through their \textit{location} in the hyperellipsoid, i.e.; they are all confined to the positive orthant. (c)~The lengths of each model manifold along the eleven axes of the hyperellipsoid \smash{$H_P$} in Eq.~\eqref{eq:vand_bnds}. Black points are the numerically computed lengths of \smash{$H_P$}, given by \smash{$2C \sqrt{N} \sigma_j(VD)$} in Eq.~\eqref{eq:vand_bnds}, and include the error term from Eq.~\eqref{eq:taylor_error} (note the kink at the second to last point), forming an upper bound on possible lengths of the manifolds. The explicit decay rate of the Chebyshev-based bound (black dotted line) is based on the fact that models obeying Eq.~\eqref{eq:deriv_bnd} are analytic in the ellipse \smash{$E_\rho(\zeta)$}. (Here, $\rho(\zeta) \approx 3.81$.) It captures the decay rate of $\sigma_j(VD)$ for $j < 11$, and closely follows the true decay rate in the successive widths of the various manifolds. } 
\label{fig:exponentials}
\end{figure*}

 
We illustrate this idea by considering three simple yet quite disparate analytic models:
\begin{enumerate}[leftmargin=*,noitemsep]
\item  \textit{Fitting exponentials}, such as for radioactive decay~\cite{Transtrum2015,Transtrum2011}. Here, we set 
\smash{$y_{\theta}(t) = \sum_{\alpha=0}^{10} A_\alpha\exp\left(-\lambda_\alpha t\right)$},
where model parameters are the amplitudes \smash{$A_\alpha$} and decay rates \smash{$\lambda_\alpha$}, and $t$ represents time.

\item \textit{Reaction velocities} of an enzyme-catalysed chemical reaction~\cite{Averick1992,Kowalik1968}. This model can be expressed as
\smash{$y_{\theta}(t) = (\theta_1 t^2 + \theta_2 t)/(t^2 + \theta_3 t + \theta_4)$}~\cite{Transtrum2015},
where $t$ represents the substrate concentration. 

\item \textit{The infected fraction of a population} in an SIR epidemiology model~\cite{Hethcote200}. This model predicts the size of a populations that is susceptible to infection ($S(t)$), infected ($I(t)$), and recovered from infection ($R(t)$). These are expressed through three coupled differential equations: \smash{$ \dot{S} = -\beta I S / N_{tot},$} \smash{$\dot{I} = \beta I S / N_{tot} - \gamma I$}, and \smash{$\dot{R} = \gamma I$},
where model parameters $\beta$ and $\gamma$ represent the rates of infection and recovery, and additional parameters include the total population $N_{tot}$, and initial infected and recovered population. At all times, $S(t)$, $I(t)$ and $R(t)$ sum to $N_{tot}$, and we set \smash{$y_{\theta}(t)=I(t)$}. 
\end{enumerate}
Fig.~\ref{fig:exponentials} displays the model manifolds for these three models. Two-dimensional projections of the manifolds are shown to be bounded within the Taylor-based hyperellipsoid \smash{$H_Y$}. The hyperribbon structure of the manifolds is accurately captured by the numerical bound from Eq.~\eqref{eq:vand_bnds}, and the decay in the bounds are clearly captured by the Chebyshev rate from Eq.~\eqref{eq:HY}.  

\vspace{.25cm}

\textbf{ 2D models.} Experiments are often conducted with more than one experimental condition, such as time and temperature. Consider the 2D model \smash{$y_\theta(t,s)$}, where model predictions \smash{$Y_{jk}$} are made at points \smash{$t_j, s_k$}, with $0 \leq j\!+\!k\! \leq N\!-\!1$. We can again use polynomial  approximation to constrain the geometry of the resulting model manifold $\mathcal{Y}$. 
In this case, we assume without loss of generality that \smash{$(t, s) \in [-1, 1]^2$}, and we assume \smash{$y_\theta$} can be expressed as a 2D Chebyshev expansion:
\smash{$y_\theta (t, s) = \sum_{j = 0}^\infty \sum_{k = 0}^\infty c_{jk}(\theta) T_{jk}(t,s)$}, where \smash{$T_{jk}(t, s) = T_j(t)T_k(s)$}. 
The following 2D polynomial of total degree $N\!-\!1$ approximates \smash{$y_\theta$}:
\begin{equation} 
 p_{N-1}(t, s; \theta) = \sum_{0 \leq j+k \leq N-1}c_{jk}(\theta) T_{jk}(t,s).
\end{equation} 
Let $\rho > 1$ and $M > 0$ be constants. For all fixed choices of \smash{$s=s^*$}, suppose that the 1D function of $t$,  \smash{$y_\theta(t, s^*)$}, is analytic in $t$ and bounded $\leq M$ uniformly with respect to both $s$ and $\theta$, and that an analogous condition holds for \smash{$y_\theta(s, t^*)$}. A result similar to Theorem~\ref{thm:Analytic_conv} can be proven by adapting the ideas in~\cite[Ch.~8]{trefethen2013approximation} to the 2D setting. Specifically, we have that 
\begin{align} 
(i) \, \,  & \|y-p_{N-1}\|_\infty \leq 4MN C_1\rho^{-N+1}, \label{eq:2Derror} \\ 
(ii) \, \, & |c_{jk}(\theta)| \leq 4M \rho^{-(j+k)}, \label{eq:2Dcoeffs}
\end{align} 
where \smash{$C_1= (2\rho-1)/(1-\rho)^2$}.

As in the 1D case, we study the model manifold $\mathcal{P}$ associated with \smash{$p_{N-1}$} as an approximation to $\mathcal{Y}$, the manifold for \smash{$y_\theta$}.  We parameterize $\mathcal{P}$ using a vector of blocks,  \smash{$P(\theta)  = (B_0,\ldots, B_{N-1})^T$}, where \smash{$B_j = (P_{0j}, P_{1(j-1)}, \ldots, P_{j0})$} and \smash{$P_{jk} = p_{N-1}(t_j, s_k; \theta)$}.  Since each block \smash{$B_j$} has $j\!+\!1$ entries, $P(\theta)$ is of length $n = N(N\!+\!1)/2$. Corresponding vectors of sample locations $\vec{t}$ and $\vec{s}$ are defined so that   \smash{$P(\theta) = p_{N-1}(\vec{t}, \vec{s}; \theta)$}. 

As before, we exploit the decay of the bounds in Eq.~\eqref{eq:2Dcoeffs} to show that  $P$ lies in the range of a matrix with strongly decaying singular values. To see this, define \smash{$\vec{\tilde{c}}$} as an appropriately ordered $n \times 1$ vector of the scaled coefficients \smash{$\tilde{c}_{jk} = \rho^{-(j+k)}c_{jk}$}, and form the linear map \smash{$P(\theta) = X \vec{\tilde{c}}$}. Here, \smash{$ X= [ X_{B^0} | \cdots | X_{B^{N-1}}]$}, where \smash{$X_{B^j}$} is a block of $j\!+\!1$ columns scaled by \smash{$\rho^{-j}$}. Specifically, 
\smash{$X_{B^j} = \rho^{-j} [T_0(\vec{t})T_{j}(\vec{s})\, | T_{1}(\vec{t})T_{j-1}(\vec{s})\, |\cdots \,| T_{j}(\vec{t}) T_0(\vec{s})]$}.   Since \smash{$\vec{\tilde{c}}$}  is constrained to lie in an $n$-sphere of radius \smash{$4M\sqrt{n}$}, the manifold $\mathcal{P}$ is contained in a hyperellipsoid \smash{$H_P$}  with cross-sectional widths characterized by the singular values of $X$.
One can show that the singular values of $X$ must decay at, at least, a subgeometric rate. 
An argument similar to the one used in Theorem~\ref{thm: ESE} shows that for $2 \leq j \leq n$, 
\begin{equation}
\label{eq:2Dsvs}
 \sigma_j(X) \leq \frac{3\sqrt{C_2} }{2}n \rho^{-\left \lfloor  \sqrt{8(j-1)+1} /2  -1/2\right \rfloor},
\end{equation} 
where \smash{$C_2= (1 + \rho^{-2} + \rho^{-4})/(1-\rho^{-2})^3$} and \smash{$\lfloor \, \cdot \,  \rfloor$} represents the floor function.  One can use \smash{$H_P$} and Eq.~\eqref{eq:2Derror} to explicitly construct a hyperellipsoid \smash{$H_Y$} that must contain $\mathcal{Y}$, just as in Eq.~\eqref{eq:dim_bound2}.  
 We expand our three previous models to the 2D setting in the supplementary materials to illustrate our bounds. While our results are stated in terms of Chebyshev expansions, a similar argument can be made using 2D Taylor expansions, and all of these ideas extend naturally to the multidimensional case.

A fundamental characteristic of the universality class of sloppy models is their global geometry: their model manifolds are shaped like hyperribbons. Through polynomial approximation, we obtained a bound on the size and shape of their manifold that accurately captures the relative widths of these hyperribbons. This bound, a hyperellipsoid, is controlled by the analyticity and smoothness of the underlying models.  Our results therefore establish a rigorous framework that explains the role of model smoothness in the observation of sloppiness.


\textbf{Acknowledgements.} We thank Mark Transtrum for suggestions related to selecting models used in this letter. KNQ was supported by a fellowship from the Natural Sciences and Engineering Research Council of Canada (NSERC), and JPS and KNQ were supported by the National Science Foundation (NSF) through grant DMR-1719490. AT was supported by NSF grant no.~DMS-1818757, and HW was supported by NSF grant no.~DGE-1650441.

\bibliography{Chebyshev_SloppyModels}

\begin{thebibliography}{20}%
\makeatletter
\providecommand \@ifxundefined [1]{%
 \@ifx{#1\undefined}
}%
\providecommand \@ifnum [1]{%
 \ifnum #1\expandafter \@firstoftwo
 \else \expandafter \@secondoftwo
 \fi
}%
\providecommand \@ifx [1]{%
 \ifx #1\expandafter \@firstoftwo
 \else \expandafter \@secondoftwo
 \fi
}%
\providecommand \natexlab [1]{#1}%
\providecommand \enquote  [1]{``#1''}%
\providecommand \bibnamefont  [1]{#1}%
\providecommand \bibfnamefont [1]{#1}%
\providecommand \citenamefont [1]{#1}%
\providecommand \href@noop [0]{\@secondoftwo}%
\providecommand \href [0]{\begingroup \@sanitize@url \@href}%
\providecommand \@href[1]{\@@startlink{#1}\@@href}%
\providecommand \@@href[1]{\endgroup#1\@@endlink}%
\providecommand \@sanitize@url [0]{\catcode `\\12\catcode `\$12\catcode
  `\&12\catcode `\#12\catcode `\^12\catcode `\_12\catcode `\%12\relax}%
\providecommand \@@startlink[1]{}%
\providecommand \@@endlink[0]{}%
\providecommand \url  [0]{\begingroup\@sanitize@url \@url }%
\providecommand \@url [1]{\endgroup\@href {#1}{\urlprefix }}%
\providecommand \urlprefix  [0]{URL }%
\providecommand \Eprint [0]{\href }%
\providecommand \doibase [0]{http://dx.doi.org/}%
\providecommand \selectlanguage [0]{\@gobble}%
\providecommand \bibinfo  [0]{\@secondoftwo}%
\providecommand \bibfield  [0]{\@secondoftwo}%
\providecommand \translation [1]{[#1]}%
\providecommand \BibitemOpen [0]{}%
\providecommand \bibitemStop [0]{}%
\providecommand \bibitemNoStop [0]{.\EOS\space}%
\providecommand \EOS [0]{\spacefactor3000\relax}%
\providecommand \BibitemShut  [1]{\csname bibitem#1\endcsname}%
\let\auto@bib@innerbib\@empty
\bibitem [{\citenamefont {Brown}\ and\ \citenamefont
  {Sethna}(2003)}]{Transtrum2003}%
  \BibitemOpen
  \bibfield  {author} {\bibinfo {author} {\bibfnamefont {K.}~\bibnamefont
  {Brown}}\ and\ \bibinfo {author} {\bibfnamefont {J.}~\bibnamefont {Sethna}},\
  }\href@noop {} {\bibfield  {journal} {\bibinfo  {journal} {Phys. Rev. E}\
  }\textbf {\bibinfo {volume} {68}} (\bibinfo {year} {2003})}\BibitemShut
  {NoStop}%
\bibitem [{\citenamefont {Brown}\ \emph {et~al.}(2004)\citenamefont {Brown},
  \citenamefont {Hill}, \citenamefont {Calero}, \citenamefont {Myers},
  \citenamefont {Lee}, \citenamefont {Sethna},\ and\ \citenamefont
  {Cerione}}]{Brown2004}%
  \BibitemOpen
  \bibfield  {author} {\bibinfo {author} {\bibfnamefont {K.}~\bibnamefont
  {Brown}}, \bibinfo {author} {\bibfnamefont {C.}~\bibnamefont {Hill}},
  \bibinfo {author} {\bibfnamefont {C.}~\bibnamefont {Calero}}, \bibinfo
  {author} {\bibfnamefont {C.}~\bibnamefont {Myers}}, \bibinfo {author}
  {\bibfnamefont {K.}~\bibnamefont {Lee}}, \bibinfo {author} {\bibfnamefont
  {J.}~\bibnamefont {Sethna}}, \ and\ \bibinfo {author} {\bibfnamefont {R.~A.}\
  \bibnamefont {Cerione}},\ }\href@noop {} {\bibfield  {journal} {\bibinfo
  {journal} {Phys. Biol.}\ }\textbf {\bibinfo {volume} {1}} (\bibinfo {year}
  {2004})}\BibitemShut {NoStop}%
\bibitem [{\citenamefont {Gutenkunst}\ \emph {et~al.}(2007)\citenamefont
  {Gutenkunst}, \citenamefont {Waterfall}, \citenamefont {Casey}, \citenamefont
  {Brown}, \citenamefont {Myers},\ and\ \citenamefont
  {Sethna}}]{Gutenkunst_SloppyParameters_2007}%
  \BibitemOpen
  \bibfield  {author} {\bibinfo {author} {\bibfnamefont {R.~N.}\ \bibnamefont
  {Gutenkunst}}, \bibinfo {author} {\bibfnamefont {J.~J.}\ \bibnamefont
  {Waterfall}}, \bibinfo {author} {\bibfnamefont {F.~P.}\ \bibnamefont
  {Casey}}, \bibinfo {author} {\bibfnamefont {K.~S.}\ \bibnamefont {Brown}},
  \bibinfo {author} {\bibfnamefont {C.~R.}\ \bibnamefont {Myers}}, \ and\
  \bibinfo {author} {\bibfnamefont {J.~P.}\ \bibnamefont {Sethna}},\ }\href
  {\doibase 10.1371/journal.pcbi.0030189} {\bibfield  {journal} {\bibinfo
  {journal} {PLOS Comput. Bio}\ }\textbf {\bibinfo {volume} {3}},\ \bibinfo
  {pages} {1} (\bibinfo {year} {2007})}\BibitemShut {NoStop}%
\bibitem [{\citenamefont {Gutenkunst}(2007)}]{Ryan2007}%
  \BibitemOpen
  \bibfield  {author} {\bibinfo {author} {\bibfnamefont {R.}~\bibnamefont
  {Gutenkunst}},\ }\emph {\bibinfo {title} {Sloppiness, modeling, and evolution
  in biochemical networks}},\ \href {http://ecommons.library.
  cornell.edu/handle/1813/8206} {Ph.D. thesis},\ \bibinfo  {school} {Cornell
  University} (\bibinfo {year} {2007})\BibitemShut {NoStop}%
\bibitem [{\citenamefont {Ruhe}(1980)}]{Ruhe1980}%
  \BibitemOpen
  \bibfield  {author} {\bibinfo {author} {\bibfnamefont {A.}~\bibnamefont
  {Ruhe}},\ }\href@noop {} {\bibfield  {journal} {\bibinfo  {journal} {SIAM J.
  Sci. Stat. Comput.}\ }\textbf {\bibinfo {volume} {1}} (\bibinfo {year}
  {1980})}\BibitemShut {NoStop}%
\bibitem [{\citenamefont {Machta}\ \emph {et~al.}(2013)\citenamefont {Machta},
  \citenamefont {Chachra}, \citenamefont {Transtrum},\ and\ \citenamefont
  {Sethna}}]{Machta_ParameterCompression_2013}%
  \BibitemOpen
  \bibfield  {author} {\bibinfo {author} {\bibfnamefont {B.~B.}\ \bibnamefont
  {Machta}}, \bibinfo {author} {\bibfnamefont {R.}~\bibnamefont {Chachra}},
  \bibinfo {author} {\bibfnamefont {M.~K.}\ \bibnamefont {Transtrum}}, \ and\
  \bibinfo {author} {\bibfnamefont {J.~P.}\ \bibnamefont {Sethna}},\ }\href
  {\doibase 10.1126/science.1238723} {\bibfield  {journal} {\bibinfo  {journal}
  {Science}\ }\textbf {\bibinfo {volume} {342}},\ \bibinfo {pages} {604}
  (\bibinfo {year} {2013})}\BibitemShut {NoStop}%
\bibitem [{\citenamefont {Berman}\ and\ \citenamefont
  {Wang}(2007)}]{Berman2007}%
  \BibitemOpen
  \bibfield  {author} {\bibinfo {author} {\bibfnamefont {G.}~\bibnamefont
  {Berman}}\ and\ \bibinfo {author} {\bibfnamefont {Z.}~\bibnamefont {Wang}},\
  }\href@noop {} {\bibfield  {journal} {\bibinfo  {journal} {J. Fluid Mech.}\
  }\textbf {\bibinfo {volume} {582}} (\bibinfo {year} {2007})}\BibitemShut
  {NoStop}%
\bibitem [{\citenamefont {Transtrum}\ \emph {et~al.}(2015)\citenamefont
  {Transtrum}, \citenamefont {Machta}, \citenamefont {Brown}, \citenamefont
  {Daniels}, \citenamefont {Myers},\ and\ \citenamefont
  {Sethna}}]{Transtrum2015}%
  \BibitemOpen
  \bibfield  {author} {\bibinfo {author} {\bibfnamefont {M.~K.}\ \bibnamefont
  {Transtrum}}, \bibinfo {author} {\bibfnamefont {B.~B.}\ \bibnamefont
  {Machta}}, \bibinfo {author} {\bibfnamefont {K.~S.}\ \bibnamefont {Brown}},
  \bibinfo {author} {\bibfnamefont {B.~C.}\ \bibnamefont {Daniels}}, \bibinfo
  {author} {\bibfnamefont {C.~R.}\ \bibnamefont {Myers}}, \ and\ \bibinfo
  {author} {\bibfnamefont {J.~P.}\ \bibnamefont {Sethna}},\ }\href {\doibase
  10.1063/1.4923066} {\bibfield  {journal} {\bibinfo  {journal} {J. Chem.
  Phys}\ }\textbf {\bibinfo {volume} {143}} (\bibinfo {year} {2015}),\
  10.1063/1.4923066},\ \Eprint {http://arxiv.org/abs/1501.07668} {1501.07668}
  \BibitemShut {NoStop}%
\bibitem [{\citenamefont {Trefethen}(2013)}]{trefethen2013approximation}%
  \BibitemOpen
  \bibfield  {author} {\bibinfo {author} {\bibfnamefont {L.~N.}\ \bibnamefont
  {Trefethen}},\ }\href@noop {} {\emph {\bibinfo {title} {Approximation
  {T}heory and {A}pproximation {P}ractice}}}\ (\bibinfo  {publisher} {SIAM},\
  \bibinfo {year} {2013})\BibitemShut {NoStop}%
\bibitem [{\citenamefont {Transtrum}\ and\ \citenamefont
  {Qiu}(2014)}]{Transtrum2014}%
  \BibitemOpen
  \bibfield  {author} {\bibinfo {author} {\bibfnamefont {M.~K.}\ \bibnamefont
  {Transtrum}}\ and\ \bibinfo {author} {\bibfnamefont {P.}~\bibnamefont
  {Qiu}},\ }\href {\doibase 10.1103/PhysRevLett.113.098701} {\bibfield
  {journal} {\bibinfo  {journal} {PRL}\ }\textbf {\bibinfo {volume} {113}},\
  \bibinfo {pages} {1} (\bibinfo {year} {2014})}\BibitemShut {NoStop}%
\bibitem [{\citenamefont {Transtrum}\ \emph {et~al.}(2011)\citenamefont
  {Transtrum}, \citenamefont {Machta},\ and\ \citenamefont
  {Sethna}}]{Transtrum2011}%
  \BibitemOpen
  \bibfield  {author} {\bibinfo {author} {\bibfnamefont {M.}~\bibnamefont
  {Transtrum}}, \bibinfo {author} {\bibfnamefont {B.~B.}\ \bibnamefont
  {Machta}}, \ and\ \bibinfo {author} {\bibfnamefont {J.}~\bibnamefont
  {Sethna}},\ }\href@noop {} {\bibfield  {journal} {\bibinfo  {journal} {Phys.
  Rev. E}\ }\textbf {\bibinfo {volume} {83}} (\bibinfo {year}
  {2011})}\BibitemShut {NoStop}%
\bibitem [{\citenamefont {Transtrum}\ \emph {et~al.}(2010)\citenamefont
  {Transtrum}, \citenamefont {Machta},\ and\ \citenamefont
  {Sethna}}]{Transtrum2010}%
  \BibitemOpen
  \bibfield  {author} {\bibinfo {author} {\bibfnamefont {M.~K.}\ \bibnamefont
  {Transtrum}}, \bibinfo {author} {\bibfnamefont {B.~B.}\ \bibnamefont
  {Machta}}, \ and\ \bibinfo {author} {\bibfnamefont {J.~P.}\ \bibnamefont
  {Sethna}},\ }\href@noop {} {\bibfield  {journal} {\bibinfo  {journal} {Phys.
  Rev. Lett.}\ }\textbf {\bibinfo {volume} {104}} (\bibinfo {year}
  {2010})}\BibitemShut {NoStop}%
\bibitem [{\citenamefont {Waterfall}\ \emph {et~al.}(2006)\citenamefont
  {Waterfall}, \citenamefont {Casey}, \citenamefont {Gutenkunst}, \citenamefont
  {Brown}, \citenamefont {Myers}, \citenamefont {Brouwer}, \citenamefont
  {Elser},\ and\ \citenamefont {Sethna}}]{Waterfall2006}%
  \BibitemOpen
  \bibfield  {author} {\bibinfo {author} {\bibfnamefont {J.~J.}\ \bibnamefont
  {Waterfall}}, \bibinfo {author} {\bibfnamefont {F.~P.}\ \bibnamefont
  {Casey}}, \bibinfo {author} {\bibfnamefont {R.~N.}\ \bibnamefont
  {Gutenkunst}}, \bibinfo {author} {\bibfnamefont {K.~S.}\ \bibnamefont
  {Brown}}, \bibinfo {author} {\bibfnamefont {C.~R.}\ \bibnamefont {Myers}},
  \bibinfo {author} {\bibfnamefont {P.~W.}\ \bibnamefont {Brouwer}}, \bibinfo
  {author} {\bibfnamefont {V.}~\bibnamefont {Elser}}, \ and\ \bibinfo {author}
  {\bibfnamefont {J.~P.}\ \bibnamefont {Sethna}},\ }\href {\doibase
  10.1103/PhysRevLett.97.150601} {\bibfield  {journal} {\bibinfo  {journal}
  {PRL}\ }\textbf {\bibinfo {volume} {97}},\ \bibinfo {pages} {150601}
  (\bibinfo {year} {2006})}\BibitemShut {NoStop}%
\bibitem [{Note1()}]{Note1}%
  \BibitemOpen
  \bibinfo {note} {Previous proofs with weaker bounds were provided through
  private communications with Ari Turner and Yaming Yu.}\BibitemShut {Stop}%
\bibitem [{\citenamefont {H.~Golub}\ and\ \citenamefont
  {Van~Loan}(1996)}]{Gene1996}%
  \BibitemOpen
  \bibfield  {author} {\bibinfo {author} {\bibfnamefont {G.}~\bibnamefont
  {H.~Golub}}\ and\ \bibinfo {author} {\bibfnamefont {C.~F.}\ \bibnamefont
  {Van~Loan}},\ }\href@noop {} {\emph {\bibinfo {title} {Matrix
  Computations}}}\ (\bibinfo  {publisher} {Johns Hopkins University Press,
  Baltimore},\ \bibinfo {year} {1996})\BibitemShut {NoStop}%
\bibitem [{\citenamefont {Beckermann}\ and\ \citenamefont
  {Townsend}(2017)}]{beckermann2017singular}%
  \BibitemOpen
  \bibfield  {author} {\bibinfo {author} {\bibfnamefont {B.}~\bibnamefont
  {Beckermann}}\ and\ \bibinfo {author} {\bibfnamefont {A.}~\bibnamefont
  {Townsend}},\ }\href@noop {} {\bibfield  {journal} {\bibinfo  {journal} {SIAM
  J. Matrix Anal. \& Appl.}\ }\textbf {\bibinfo {volume} {38}},\ \bibinfo
  {pages} {1227} (\bibinfo {year} {2017})}\BibitemShut {NoStop}%
\bibitem [{\citenamefont {Demanet}\ and\ \citenamefont
  {Townsend}(2018)}]{demanet2016stable}%
  \BibitemOpen
  \bibfield  {author} {\bibinfo {author} {\bibfnamefont {L.}~\bibnamefont
  {Demanet}}\ and\ \bibinfo {author} {\bibfnamefont {A.}~\bibnamefont
  {Townsend}},\ }\href {\doibase 10.1007/s10208-018-9384-1} {\bibfield
  {journal} {\bibinfo  {journal} {Found. Comput. Math.}\ } (\bibinfo {year}
  {2018}),\ 10.1007/s10208-018-9384-1}\BibitemShut {NoStop}%
\bibitem [{\citenamefont {Averick}\ \emph {et~al.}(1992)\citenamefont
  {Averick}, \citenamefont {Carter},\ and\ \citenamefont {Xue}}]{Averick1992}%
  \BibitemOpen
  \bibfield  {author} {\bibinfo {author} {\bibfnamefont {B.}~\bibnamefont
  {Averick}}, \bibinfo {author} {\bibfnamefont {J.~M.}\ \bibnamefont {Carter}},
  \ and\ \bibinfo {author} {\bibfnamefont {G.}~\bibnamefont {Xue}},\
  }\href@noop {} {\bibfield  {journal} {\bibinfo  {journal} {Preprint
  MCS-P153-0694, Mathematics and Computer Science Division, Argonne National
  Laboratory, Argonne, Illinois}\ } (\bibinfo {year} {1992})}\BibitemShut
  {NoStop}%
\bibitem [{\citenamefont {Kowalik}\ and\ \citenamefont
  {Morrison}(1968)}]{Kowalik1968}%
  \BibitemOpen
  \bibfield  {author} {\bibinfo {author} {\bibfnamefont {J.}~\bibnamefont
  {Kowalik}}\ and\ \bibinfo {author} {\bibfnamefont {J.}~\bibnamefont
  {Morrison}},\ }\href@noop {} {\bibfield  {journal} {\bibinfo  {journal}
  {Math. Biosci.}\ }\textbf {\bibinfo {volume} {2}} (\bibinfo {year}
  {1968})}\BibitemShut {NoStop}%
\bibitem [{\citenamefont {Hethcote}(2000)}]{Hethcote200}%
  \BibitemOpen
  \bibfield  {author} {\bibinfo {author} {\bibfnamefont {H.~W.}\ \bibnamefont
  {Hethcote}},\ }\href@noop {} {\bibfield  {journal} {\bibinfo  {journal} {SIAM
  Review}\ }\textbf {\bibinfo {volume} {42}},\ \bibinfo {pages} {599} (\bibinfo
  {year} {2000})}\BibitemShut {NoStop}%
\end{thebibliography}%


\begin{thebibliography}{3}%
\makeatletter
\providecommand \@ifxundefined [1]{%
 \@ifx{#1\undefined}
}%
\providecommand \@ifnum [1]{%
 \ifnum #1\expandafter \@firstoftwo
 \else \expandafter \@secondoftwo
 \fi
}%
\providecommand \@ifx [1]{%
 \ifx #1\expandafter \@firstoftwo
 \else \expandafter \@secondoftwo
 \fi
}%
\providecommand \natexlab [1]{#1}%
\providecommand \enquote  [1]{``#1''}%
\providecommand \bibnamefont  [1]{#1}%
\providecommand \bibfnamefont [1]{#1}%
\providecommand \citenamefont [1]{#1}%
\providecommand \href@noop [0]{\@secondoftwo}%
\providecommand \href [0]{\begingroup \@sanitize@url \@href}%
\providecommand \@href[1]{\@@startlink{#1}\@@href}%
\providecommand \@@href[1]{\endgroup#1\@@endlink}%
\providecommand \@sanitize@url [0]{\catcode `\\12\catcode `\$12\catcode
  `\&12\catcode `\#12\catcode `\^12\catcode `\_12\catcode `\%12\relax}%
\providecommand \@@startlink[1]{}%
\providecommand \@@endlink[0]{}%
\providecommand \url  [0]{\begingroup\@sanitize@url \@url }%
\providecommand \@url [1]{\endgroup\@href {#1}{\urlprefix }}%
\providecommand \urlprefix  [0]{URL }%
\providecommand \Eprint [0]{\href }%
\providecommand \doibase [0]{http://dx.doi.org/}%
\providecommand \selectlanguage [0]{\@gobble}%
\providecommand \bibinfo  [0]{\@secondoftwo}%
\providecommand \bibfield  [0]{\@secondoftwo}%
\providecommand \translation [1]{[#1]}%
\providecommand \BibitemOpen [0]{}%
\providecommand \bibitemStop [0]{}%
\providecommand \bibitemNoStop [0]{.\EOS\space}%
\providecommand \EOS [0]{\spacefactor3000\relax}%
\providecommand \BibitemShut  [1]{\csname bibitem#1\endcsname}%
\let\auto@bib@innerbib\@empty
\bibitem [{\citenamefont {Trefethen}(2013)}]{trefethen2013approximation}%
  \BibitemOpen
  \bibfield  {author} {\bibinfo {author} {\bibfnamefont {L.~N.}\ \bibnamefont
  {Trefethen}},\ }\href@noop {} {\emph {\bibinfo {title} {Approximation
  {T}heory and {A}pproximation {P}ractice}}}\ (\bibinfo  {publisher} {SIAM},\
  \bibinfo {year} {2013})\BibitemShut {NoStop}%
\bibitem [{\citenamefont {Weideman}\ and\ \citenamefont
  {Trefethen}(2007)}]{weideman2007kink}%
  \BibitemOpen
  \bibfield  {author} {\bibinfo {author} {\bibfnamefont {J.}~\bibnamefont
  {Weideman}}\ and\ \bibinfo {author} {\bibfnamefont {L.~N.}\ \bibnamefont
  {Trefethen}},\ }\href@noop {} {\bibfield  {journal} {\bibinfo  {journal}
  {Numerische Mathematik}\ }\textbf {\bibinfo {volume} {107}},\ \bibinfo
  {pages} {707} (\bibinfo {year} {2007})}\BibitemShut {NoStop}%
\bibitem [{\citenamefont {Townsend}\ and\ \citenamefont
  {Wilber}(2018)}]{townsend2018singular}%
  \BibitemOpen
  \bibfield  {author} {\bibinfo {author} {\bibfnamefont {A.}~\bibnamefont
  {Townsend}}\ and\ \bibinfo {author} {\bibfnamefont {H.}~\bibnamefont
  {Wilber}},\ }\href@noop {} {\bibfield  {journal} {\bibinfo  {journal} {Lin.
  Alg. \& Appl.}\ }\textbf {\bibinfo {volume} {548}},\ \bibinfo {pages} {19}
  (\bibinfo {year} {2018})}\BibitemShut {NoStop}%
\end{thebibliography}%

\end{document}


\title{Supplemental Material for: Chebyshev approximation and the global geometry of sloppy models}

\author{Katherine N. Quinn}
\affiliation{Physics Department, Cornell University, Ithaca, NY 14853-2501, United States. (knq2@cornell.edu)}

\author{Heather Wilber}
\affiliation{Center for Applied Mathematics, Cornell University, Ithaca, NY 14853-3801, United States. (hdw27@cornell.edu)}

\author{Alex Townsend}
\affiliation{Mathematics Department, Cornell University, Ithaca, NY 14853-4201, United States. (townsend@cornell.edu)}

\author{James P. Sethna}
\affiliation{Physics Department, Cornell University, Ithaca, NY 14853-2501, United States. (sethna@lassp.cornell.edu)}

\date{\today}

\maketitle

In Section~\ref{sec:NonA}, we discuss how our bounds can be extended to models that are not analytically continuable,  but are $k$-times continuously differentiable on the interval of approximation. 
In Section~\ref{sec:NumTests},  we give numerical results for high-dimensional manifolds and discuss the behavior of the singular values of the matrix $VD$ associated with truncated Taylor expansions. 
In Section~\ref{sec:2D}, we  extend the 1D models described in the main text to include two experimental conditions, and show that their manifolds exhibit a hyperribbon-like structure that is captured by our bounds. 
Finally, in Section~\ref{sec:Code}, we show how the visualizations of the model manifolds in the main text were generated.

\section{Non-Analytic Models}
\label{sec:NonA}

In the main text, we considered models \smash{$y_\theta(t)$}, \smash{$t\in [-1, 1]$}, that are continuously dependent on parameters $\theta = (\theta_1, \ldots, \theta_K)$ and analytic in an open neighborhood of $[-1, 1]$. We bounded the model manifold $\mathcal{Y}$ of model predictions by considering the truncated Chebyshev approximation
\begin{equation}
\label{eq:cheb_appx}
p_{N-1}(t; \theta) = \sum_{j = 0}^{N-1} c_j(\theta) T_j(t), 
\end{equation}
where \smash{$T_j$} is the Chebyshev polynomial of degree $j$.  
When \smash{$y_\theta$} is not analytic on $[-1, 1]$, the convergence of Eq.~\eqref{eq:cheb_appx} to \smash{$y_\theta$} as  $N \to \infty$ is still controlled by the smoothness of \smash{$y_\theta$}.  A standard result supplied in~\cite[Ch.~7]{trefethen2013approximation} states that if  \smash{$y_\theta$} has \smash{$\nu-1\geq 0$} derivatives that are absolutely continuous on $[-1, 1]$, with the $\nu$th derivative of total bounded variation $V<\infty$, then 
\begin{align*} 
&  (i) \, \,   \|y_\theta - p_{N-1}\|_\infty \leq \frac{2 V}{\pi \nu} (N-1-\nu)^{-\nu}, \quad N > \nu+1, \\ 
&  (ii) \, \, |c_j| \leq \frac{2V}{\pi} (j - \nu)^{-(\nu+1)}, \qquad j \geq \nu + 1.
\end{align*}

\begin{figure}[h!] 
 \centering
 \begin{minipage}{.40\textwidth} 
 \centering
  \begin{overpic}[width=\textwidth]{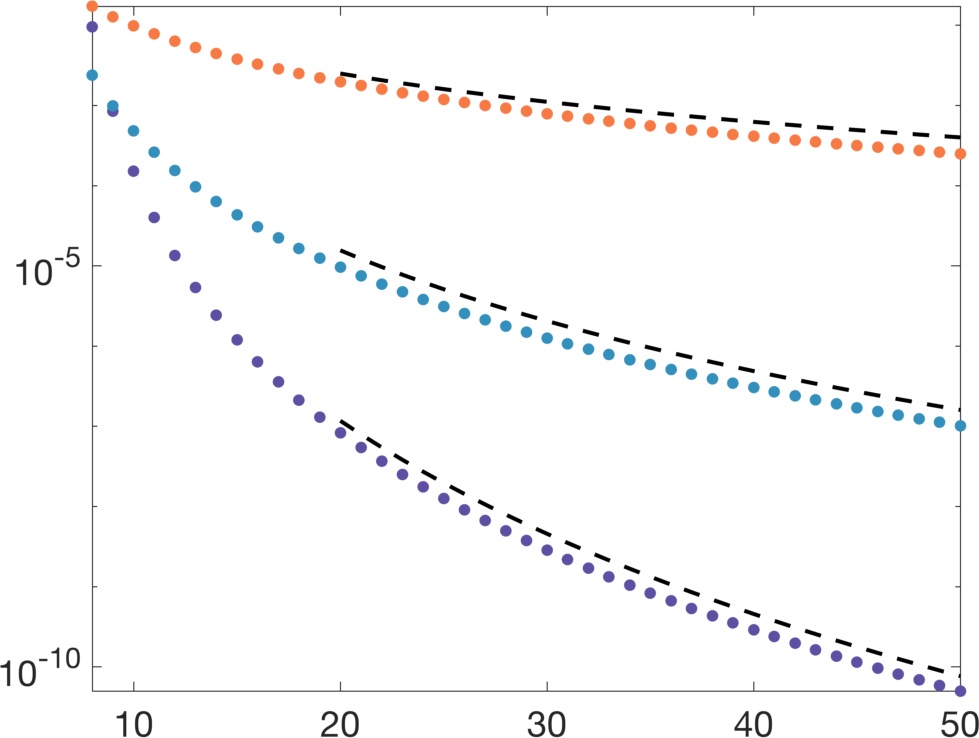}
  \put (-5, 18){\rotatebox{90}{ \small{Size of Singular Value}}}
  \put (35, -6){{ \small{Ordered Singular Value}}}
  \put (55, 23){\rotatebox{-19}{ \small{$\nu = 5, \, \mathcal{O}(j^{-8})$ }}}
  \put (60, 43){\rotatebox{-10}{ \small{$\nu = 3,  \, \mathcal{O}(j^{-5})$ }}}
  \put (68, 67){\rotatebox{-5}{ \small{$\nu = 1, \,  \mathcal{O}(j^{-2})$ }}}
  \end{overpic}
  \vspace{0.1in}
  \end{minipage}
  \caption{\textbf{The singular values \smash{$\sigma_j(X)$}}, where $X$ is described in Section~\ref{sec:NonA}, are plotted on a log scale against the index $j$ for three models of the form \smash{$y_\theta(t) = f(\theta) |t|^\nu $}: $\nu = 1$ (orange), $\nu = 3$ (blue), and $\nu= 5$ (purple). For simplicity,  we assume $f$ is smooth and independent of $t$.   In each case, the model \smash{$y_\theta$} is $\nu$-times differentiable on $[-1, 1]$. The asymptotic decay of the singular values (dotted black lines) is algebraic, with stronger decay rates as $\nu$ becomes larger. This suggests that continuously differentiable models have manifolds with (fat) hyperribbon structures, since a $\nu$-times differentiable  model \smash{$y_\theta$} has a manifold enclosed in \smash{$H_Y$}, with \smash{$\ell_j(H_Y) \approx 2 r \sigma_j(X)$} for some constant $r>0$.  }
  \label{fig:alg_decay}
 \end{figure}

To bound $\mathcal{P}$, the model manifold of \smash{$p_{N-1}(\vec{t})$}, we note that  \smash{$p_{N-1}(\vec{t}) = X  \vec{\tilde{c}}$} for \smash{$\vec{t} = (t_0, \ldots, t_{N-1})^T$}, where  \smash{$X = JD$}, with \smash{$J_{ij} = T_{j-1}(t_{i-1})$}, \smash{$D_{jj }= (j-1-\nu)^{-(\nu+1)}$} for \smash{$j \geq \nu+2$}, with \smash{$D_{jj }= 1$} otherwise. Likewise, we set \smash{$\vec{\tilde{c}} = (\tilde{c}_0, \ldots, \tilde{c}_{N-1})^T$}, where \smash{$\tilde{c}_j = (j-\nu)^{(\nu+1)}c_j$} for $j \geq \nu+1$,  and \smash{$\tilde{c}_j = c_j$} otherwise. The singular values of $X$ decay at, at least, an algebraic rate that increases with $\nu$ (see Fig.~\ref{fig:alg_decay}).  As in the analytic case, one can use $X$ as a linear map and construct a  hyperellipsoid \smash{$H_Y$} that bounds the model manifold associated with \smash{$y_\theta(\vec{t})$}.  Its cross sections are controlled by the singular values of $X$ and typically shrink algebraically fast. 

As a question of nomenclature, we suggest that an object with an algebraic
decay of widths should also be described as a hyperribbon. 

\section{Numerical observations for high dimensional manifolds}
\label{sec:NumTests}

In the main text, we bounded model predictions \smash{$y_\theta(\vec{t})$} evaluated at $N$ points \smash{$\vec{t} = (t_0,\dots,t_{N-1})^T$} by approximating \smash{$y_\theta$} with its degree $\leq N\!-\!1$ truncated Taylor expansion, which we denote by \smash{$p_{N-1}(t; \theta)$}.  The manifold associated with \smash{$p_{N-1}$} is bounded within a hyperellipsoid \smash{$H_P$}. The cross-sectional diameters of \smash{$H_P$}  are defined in terms of the singular values of the column-scaled Vandermonde matrix $X = VD$, where \smash{$(VD)_{ij} = t_{i-1}^{j-1}R^{-(j-1)}$}. Specifically, we have that
 \begin{equation} 
 \label{eq:vand_bnds}
 \ell_j(H_P) = 2 C\sqrt{N} \sigma_j(VD) \leq \frac{CN}{ \sqrt{R^2 - 1}} R^{-j + 2}, 
 \end{equation}
where $C>0$,  $R>1$ come from the analyticity constraint
\begin{eqnarray}
\label{eq:SSbound}
\sum_{k=0}^{N-1} \left(\frac{R^k}{k!}\frac{d^ky_\theta(t)}{dt^k}\right)^2 < C^2 N.
\end{eqnarray}

One can conclude, as shown in the main text, that $\mathcal{Y}$, the manifold associated with \smash{$y_\theta(\vec{t})$}, is bounded in a hyperellipsoid \smash{$H_Y$} with cross-sectional widths obeying 
$$\ell_j(H_Y) \leq \ell_j(H_P) + 2\|y_\theta - p_{N-1}\|_\infty.$$

As discussed in the main text, one expects that the decay rate \smash{$\mathcal{O}(R^{-j})$} in Eq.~\eqref{eq:vand_bnds} is weak as an upper bound on the actual ordered widths of the underlying  hyperribbon $\mathcal{Y}$.  This is related to the fact that unlike truncated Chebyshev expansions, truncated Taylor polynomials do not converge to \smash{$y_\theta$} at a rate that is asymptotically optimal for polynomial approximants (see~\cite[Ch. 12--16]{trefethen2013approximation}). 

 However, we find that the singular values \smash{$\sigma_j(VD)$} behave in a surprising way: For small to moderate $j$, the magnitude of \smash{$\sigma_j(VD)$} decays  at a rate close to the limit  predicted by Chebyshev approximation: \smash{$\mathcal{O}( \rho_{\max}^{-j})$}, where \smash{$\rho_{\max} = R + \sqrt{R^2 + 1}$}.  It is only when $j$ is sufficiently large that \smash{$\sigma_j(VD)$} appears to decay at the predicted rate \smash{$\mathcal{O}(R^{-j})$}.  We do not yet fully understand why the singular values of $VD$ decay at two distinct rates, but speculate that it may be related to the kink observed in error plots for Clenshaw--Curtis quadrature on analytic functions~\cite{weideman2007kink}.

 \begin{figure} 
 \centering
 \begin{minipage}{.40\textwidth} 
 \centering
  \begin{overpic}[width=\textwidth]{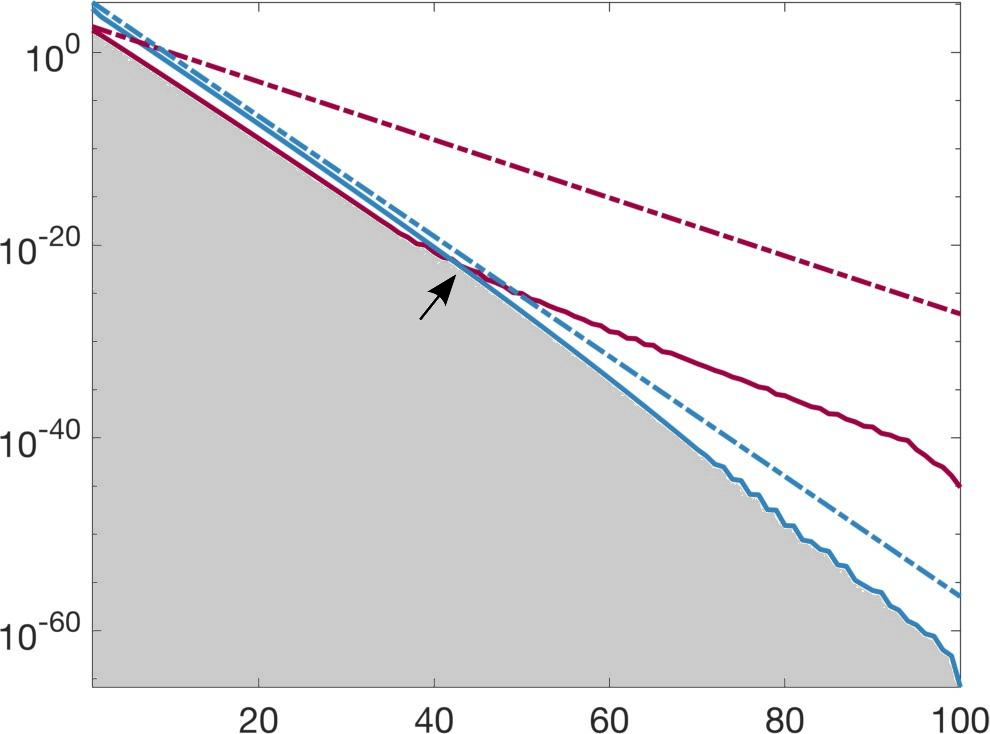}
  \put (-5, 30){\rotatebox{90}{ \small{Length}}}
  \put (25, -6){{ \small{Ordered Hyperellipsoid Axes}}}
  \put (54, 59){\rotatebox{-17}{ \small{$\mathcal{O}(R^{-j})$ }}}
  \put (37, 57){\rotatebox{-30}{ \small{$\mathcal{O}(\rho_{\max}^{-j})$ }}}
  \end{overpic}
  \vspace{0.1in}
  \end{minipage}
  \caption{\textbf{Bounds on the hyperellipsoid lengths} \smash{$\ell_j(H_P)$} using truncated Taylor (dotted purple) and truncated Chebyshev (dotted blue)  expansions are plotted on a log scale against the dimension index $j$. These form a universal bound on the ordered manifold widths of the prediction space for models \smash{$y_\theta$} that satisfy Eq.~\eqref{eq:vand_bnds}. In this example, \smash{$C = 1$}, \smash{$R = 2$}, \smash{$N = 100$}, and \smash{$\rho_{\max} \approx 4.2$}.   The solid lines show the actual computed hyperellipsoid cross-sectional lengths (on a log scale) \smash{$\ell_j(H_P)= 2 r \sigma_j(X)$}, where $X = VD$ for the Taylor-based bounds, and \smash{$X_{ij} = T_{j-1}(t_{i-1}) \rho_{\max}^{-(j-1)}$} for the Chebyshev-based bounds. The largest $40$ Taylor-based hyperellipsoid lengths decay at the  rate predicted by the Chebyshev-based bounds. Then, a kink occurs (indicated by a black arrow) and the lengths decay at the rate predicted by the bound in Eq.~\eqref{eq:vand_bnds}. For the smaller dimensions, the Chebyshev-based results produce tighter bounds. Model manifold lengths outside of the shaded region cannot occur.}
  \label{fig:bounds}
 \end{figure}

Due to this phenomenon, we find that using \smash{$\sigma_j(VD)$} directly results in good bounds on model prediction spaces for low dimensions (the larger axes of the hyperellipsoid \smash{$H_Y$}). At higher dimensions (shorter hyperellipsoid axes), the Taylor-based bounds become suboptimal, and it is beneficial to instead convert the constraint in~\eqref{eq:SSbound} to one involving Bernstein ellipses, and then use the Chebyshev-based bounds from Eq.~(8) in the main text.  The conversion of the constraint can result in bounds that are inflated by a large unphysical constant, but the decay rate in the new bound, close to \smash{$\mathcal{O}( \rho_{\max}^{-j})$}, is nearly double the rate \smash{$\mathcal{O}(R^{-j})$}. When viewed together, the Chebyshev-based bounds and numerical Taylor-based  bounds describe the successive lengths of the model  manifold across two regimes (low vs.~high dimension). We illustrate this observation using a high-dimensional manifold ($N = 100$) in Fig.~\ref{fig:bounds}.

\section{Two-Dimensional Extension of Model Predictions}
\label{sec:2D}

\begin{figure}
	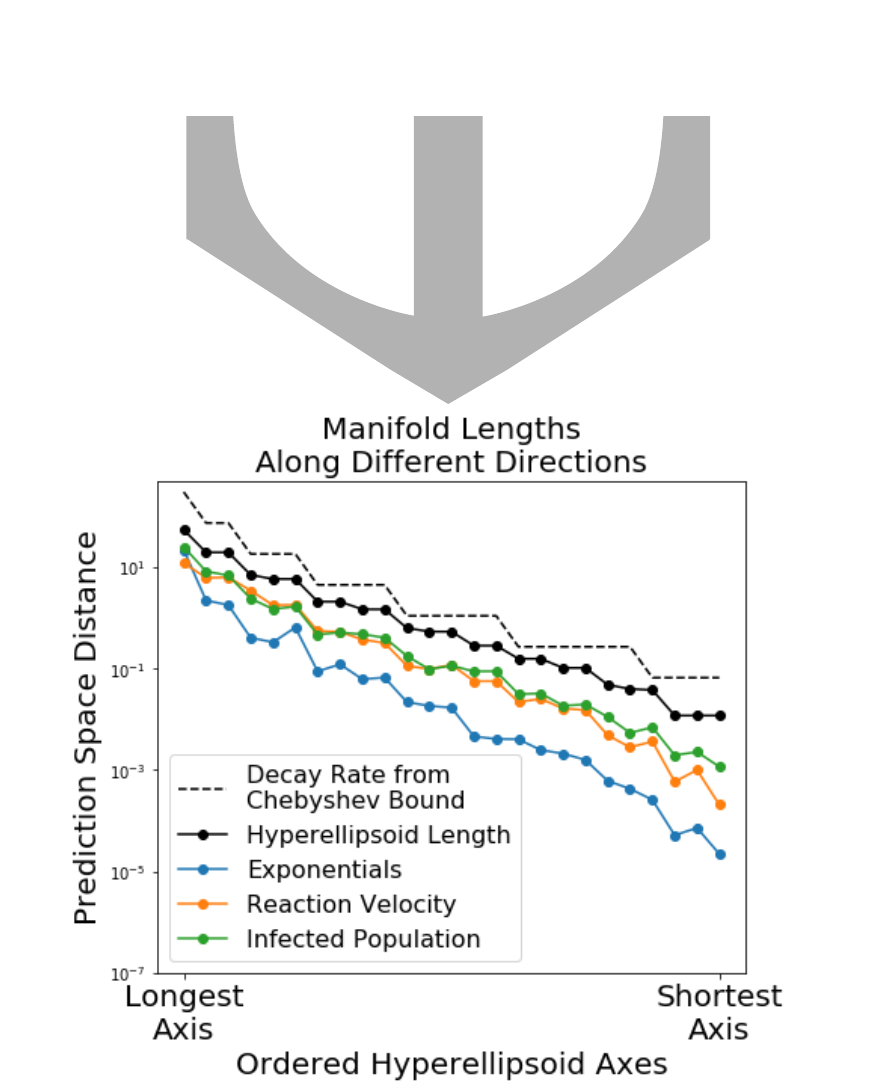
    \caption{\textbf{Model manifold} of three models with two experimental conditions: (1) exponential decay with temperature dependent decay rates, 
(2) reaction velocities of an enzyme-catalysed reaction with temperature dependent reaction rates, 
and (3) the infected population in an SIR model with infection and recovery rates that vary with parameter $s$. (a)~The models are evaluated at 25 equally spaced points \smash{$(t_i,s_i)\in [0,1]^2$} (shifted and rescaled from the interval $[-1,1]^2$) with different model parameters. All models obey the analyticity condition in Eq.~\eqref{eq:2DdiffBound} with $C$=1 and $R=2$. (b)~The explicit lengths of the three models are shown along the twenty-five axes of the hyperellipsoid \smash{$H_P$}. The upper bounds on the possible lengths (black dots) are given by \smash{$\ell_j(H_P) = 2 C \sqrt{n} \sigma_j(X)$}, where $X$ is described in Section~\ref{sec:2D}.  They exhibit subgeometric decay, with a rate that is captured by the bound in Eq.~\eqref{eq:2Dsvs} (dashed line) with $\rho \approx 4.1$.  The hierarchy of widths coming from the explicit bounds suggests that the manifolds are hyperribbons. }
    \label{fig:exponentials2D}
\end{figure}

In this section, we extend the three models used in the main text to the 2D setting. We do this by adding an extra experimental condition, denoted by $s$, to each model. In Fig.~\ref{fig:exponentials2D}, we construct the model manifolds for all three. Just as before, the model manifold is bounded by a hyperellipsoid \smash{$H_Y$} with a hierarchy of widths that form a hyperribbon structure. 
\begin{enumerate}
\item For \textit{exponentials} we consider temperature dependent decay rates,
\begin{eqnarray}
\label{eq:ExpTempDep}
\lambda_\alpha  &\rightarrow &  \lambda_\alpha \exp\left(-E_\alpha s \right), \\
y(t) &\rightarrow & y(t,s) = \sum_\alpha A_\alpha \exp\left(-\lambda_\alpha \exp (-E_\alpha s) t\right),
\end{eqnarray}
where $s = 1/T$ is inverse temperature.
\item For the model of \textit{reaction velocities}, we consider temperature dependent parameters,
\begin{eqnarray}
\label{eq:rateTempdep}
\theta_\alpha \rightarrow {\theta_\alpha} \exp\left(-E_\alpha s\right),
\end{eqnarray}
where again $s = 1/T$ is inverse temperature.
\item Finally, for the \textit{infected population} in an SIR model, we introduce infection and recovery rates that vary continuously with an infection parameter $s$ by introducing 
\begin{eqnarray}
\label{eq:SIRTempDep}
\beta &\rightarrow &  \beta \exp\left(-E_\beta s \right), \\
\gamma &\rightarrow & \gamma \exp\left(-E_\gamma s \right).
\end{eqnarray} 
\end{enumerate}
In all cases, \smash{$E_\alpha$}, \smash{$E_\beta$} and \smash{$E_\gamma$} represent activation energies in the respective models. Fig.~\ref{fig:exponentials2D} shows the model manifolds of all three example models, illustrating their hyperribbon structures. To generate these figures, we consider models that obey an analyticity constraint analogous to Eq.~\eqref{eq:SSbound}. Specifically, we assume that for all $0 \leq j\!+\!k \leq N\!-\!1$, the following condition holds uniformly in $\theta$ for a given $2D$ model \smash{$y_\theta(t,s)$}: 
\begin{eqnarray}
\label{eq:2DdiffBound}
\sum_{j+k\leq N-1}\left(\frac{R^{j+k}}{j!k!}\frac{d^{j+k}y_\B{\theta}(t,s)}{dt^{j} ds^{k}}\right)^2 < C^2 n.
\end{eqnarray}
where $R> 1, C>0$  are constants, and $n=N(N+1)/2$.  
Under this constraint, it makes sense  to bound the prediction space using truncated Taylor expansions of total degree $\leq N-1$ for small to moderate $N$ (see the discussion in Section~\ref{sec:NumTests}). This choice results in an $n \times n$ linear system of the form  \smash{$y_\theta(\vec{t}, \vec{s}) \approx  X\vec{\tilde{a}}$}, where $X$ is a column-scaled 2D Vandermonde matrix, and \smash{$\|\vec{\tilde{a}}\|_2 <C \sqrt{n}$}.  The structure of $X$ can be exploited to bound its singular values explicitly~\cite{townsend2018singular}. Alternatively, one can apply the 2D analogue to Theorem 2 from the main text to find explicit bounds in terms of $R$. In Fig.~\ref{fig:exponentials2D}, we simply use the relation  \smash{$\ell_j(H_Y) =\ell_j( H_P)+ 2\|y_\theta - p_{N-1}\|_\infty$},  and compute \smash{$\ell_j(H_P) = 2 r \sigma_j(X)$} numerically. 

We compare this with the Chebyshev-based bound established in the main text,
\begin{equation}
\label{eq:2Dsvs}
 \ell_j(H_P) \leq 2\sqrt{N}\frac{3\sqrt{C_2} }{2}n \rho^{-\left \lfloor  \sqrt{8(j-1)+1} /2  -1/2\right \rfloor},
\end{equation} 
where $\rho$ is a characteristic length related to the analyticity of the model, \smash{$C_2= (1 + \rho^{-2} + \rho^{-4})/(1-\rho^{-2})^3$}, and \smash{$\lfloor \, \cdot \,  \rfloor$} represents the floor function. This bound captures the subgeometric decay rate of the model manifold lengths for all three examples, illustrated through the dashed line in Fig.~\ref{fig:exponentials2D}.

\section{Generating model manifolds}
\label{sec:Code}

Here, we provide a detailed description of how data for the 1D models used in the main text were generated. Data for the 2D models in Section~\ref{sec:2D} were computed in a similar way.  In order to generate the model manifolds, a Monte Carlo sampling was performed on the parameter space of all three models. The model predictions for the randomly selected parameters were accepted or rejected based on whether or not they satisfied the constraint on the derivative from Eq.~\eqref{eq:SSbound}, where we set \smash{$C=1$} and \smash{$R=2$}. Since we consider eleven equally spaced points in the main manuscript, in all example models the derivative constraint was applied up to the eleventh derivative.

\begin{figure}
\vspace{.25cm}
	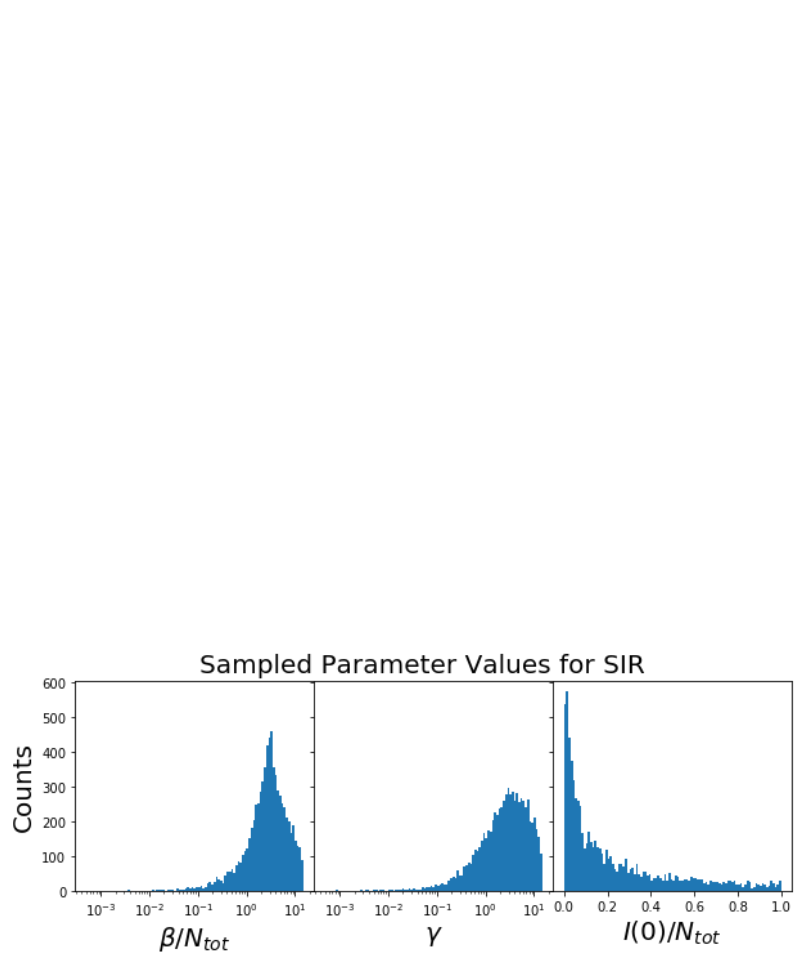
    \caption{\textbf{Histograms of valid parameter values} used to generate the model manifolds. In all the models, a Monte Carlo sampling was performed, with parameters accepted or rejected based on whether or not they satisfied the derivative condition from Eq.~\eqref{eq:SSbound}. (a)~Parameter values for exponentials, showing the distributions for the amplitudes \smash{$A_\alpha$} and decay rates \smash{$\lambda_\alpha$}. (b)~Parameter values for the reaction velocities, for each \smash{$\theta_1$, $\theta_2$, $\theta_3$} and \smash{$\theta_4$}. (c)~Parameter values for the SIR epidemiology model, showing the distribution of infection rates $\beta/N_{tot}$, recovery rates $\gamma$ and initial infected population.}
    \label{fig:paramVals}
\end{figure}

\begin{enumerate}
	\item For \textit{exponentials}, the model is of the form
	\begin{eqnarray}
	y_\theta(t) = \sum_{\alpha = 0}^{10} A_\alpha \exp\left(-\lambda_\alpha t \right),
	\end{eqnarray}
	and the derivative constraint from Eq.~\eqref{eq:SSbound} can be expressed as
	\begin{eqnarray}
	\sum_{k=0}^{N-1}\left( \sum_{\alpha=0}^{10}\frac{R^kA_\alpha}{k!}(-\lambda_\alpha)^k\exp\left(-\lambda_\alpha t\right)\right)^2  < C^2 N
	\end{eqnarray}
	for all $-1\leq t \leq 1$. From a Monte Carlo sampling, 42,000 valid samples were randomly generated. A histogram of parameters used to generate the model manifold is shown in Fig.~\ref{fig:paramVals}{(a)}.
		
	\item The model of \textit{reaction velocities} is given by
	\begin{eqnarray}
	y_{\theta}(t) = \frac{\theta_1 t^2 + \theta_2 t}{t^2 + \theta_3 t + \theta_4},
	\end{eqnarray}
	where $t$ is the substrate concentration. The derivative constraint can be expressed as
	\begin{eqnarray}
	\sum_{k=1}^N\left(\frac{R^k}{k!}\frac{d^k}{dt^k}\left(\frac{\theta_1 t^2 + \theta_2 t}{t^2 + \theta_3 t + \theta_4}\right)\right)^2 < C^2 N,
	\end{eqnarray}
	for all $-1<t<1$. We generated 24,000 valid parameter combinations, and a histogram of the different parameter values is shown in Fig.~\ref{fig:paramVals}{(b)}.
	
	\item Finally, for the \textit{infected population} in an SIR model, the number of people susceptible ($S$), infected ($I$), and recovered ($R$) are determined through three coupled differential equations:
	\begin{align*}
	&(i) \, \, \dot{S} = -\beta \frac{I S}{N_{tot}}, \\
	&(ii) \, \, \dot{I} = \beta \frac{I S}{N_{tot}} - \gamma I, \\
	&(iii) \, \, \dot{R} = \gamma I,
	\end{align*}
	where $\beta$ is the infection rate, $\gamma$ is the recovery rate, and $N_{tot}$ is the total size of the population. If we let the model predictions be the infected population, then we have \smash{$y_\theta(t) = I(t)$}. To find the $k$th derivative of such a model, we note that \smash{$\dot{S} = f_1(S,I)$} and \smash{$\dot{I} = g_1(S,I)$}. The subsequent derivatives can therefore be found recursively, by \smash{$\ddot{y}_\theta = \ddot{I}=\frac{dg_1}{dS}\dot{S} + \frac{dg_1}{dI}\dot{I}=g_2(S,I)$} and so on. From a Monte Carlo sampling, we obtained 20,000 valid parameter combinations. A histogram of parameter values used to generate the model manifold is shown in Fig.~\ref{fig:paramVals}{(c)}.
\end{enumerate}

In all three models, the smallest physically meaningful prediction is $y_\theta(t) = 0$. For exponentials and the SIR model, the largest physically meaningful prediction allowed by Eq.~\eqref{eq:SSbound} is \smash{$y_\theta(t) = C\sqrt{N}$}, and so the longest manifold distance possible is $CN$. With this sampling method, we obtained manifold lengths that are within $1.5\%$ of this maximally allowed distance, and so while more refined sampling methods could be used to resolve the manifold boundaries, they are unnecessary for our purposes.

Once a sampling of the possible parameter combinations is obtained for a model, we visualize it. Each parameter combination is evaluated at eleven equally spaced points. The space spanned by the model predictions at these points forms the model manifold $\mathcal{Y}$.

To visualize $\mathcal{Y}$, it is rotated into the basis given by the hyperellpsoid axes constructed from the space of allowed polynomials predictions, $\mathcal{P}$. Let \smash{$\{\phi_j\}_{j = 0}^{\infty}$} be a complete polynomial basis, and let \smash{$P(\vec{b})  = (P_0, \ldots, P_{N-1})$} define the model manifold $\mathcal{P}$ of \smash{$p_{N-1}(t) = \sum_{j = 0}^{N-1} b_j \phi_j(t)$}.  Polynomial predictions are given by \smash{$P_k = p_{N-1}(t_k)$}. By definition, \smash{$P(\vec{b}) = X\vec{b}$}, where  \smash{$X_{ij} = \phi_{j-1}(t_{i-1})$} and \smash{$\vec{b} = (b_0, \ldots, b_{N-1})^T$}. To find the rotation matrix used to visualize the model manifold $\mathcal{Y}$, we perform a singular value decomposition on $X$,
\begin{eqnarray}
X = U \Sigma V^T,
\end{eqnarray}
to extract the rotation matrix $U$. The data points on the model manifold are then rotated using this matrix, and visualized in Fig.~I(b) in the main text where we set $X=VD$ to be the column-scaled Vandermonde matrix.

\bibliography{Chebyshev_SloppyModels}